\documentclass[10pt, draft]{amsart}
\usepackage{amsmath}
\usepackage{latexsym}
\usepackage[all]{xy}
\numberwithin{equation}{section}

\title[Convex Invariance]{On the convex invariance in Finsler geometry}%
\author[Gallego Torrom\'e]{R. Gallego Torrom\'e }
\thanks{The author was financially supported by FAPESP, grant n. 2010/11934-6. }
\email{rgallegot@gmx.de}
\address{Departments de Matem\'atica\hfill\break\indent Universidade de S\~ao Paulo\hfill\break\indent Brazil}
\date{November 14th, 2012}

\begin{document}

\newtheorem{theorem}{Theorem}
\newtheorem{assertion}[theorem]{Assertion}
\newtheorem{proposition}[theorem]{Proposition}
\newtheorem{lemma}[theorem]{Lemma}
\newtheorem{definition}[theorem]{Definition}
\newtheorem{claim}[theorem]{Claim}
\newtheorem{corollary}[theorem]{Corollary}
\newtheorem{observation}[theorem]{Observation}
\newtheorem{conjecture}[theorem]{Conjecture}
\newtheorem{question}[theorem]{Question}
\newtheorem{example}[theorem]{Example}

\newtheorem*{theorem-A}{Theorem A}

\theoremstyle{remark}\newtheorem*{remark}{Remark}
\maketitle

\begin{abstract}
 As a natural application of the {\it theory of geometric averaging} in Finsler geometry and generalized Finsler geometry, a new approach to investigate {\it generalized Finsler geometry}, based on a convex invariance of the average structures, is introduced.
\end{abstract}

\section{Introduction}
\subsection*{Motivation}
Given a physical observable $A$, it is interesting to associate to an {\it average valued} $\langle A\rangle$, since in some cases is a reasonable expected value for the result of an experimental measurement of $A$. There are powerful reasons for this. For systems in {\it equilibrium state}, the average value of a physical quantity corresponds to the expected (most likely) value of a measurement. Also, measurement experiments are linked with a large number of identically prepared systems. Then the {\it Central Limit Theorem} (see for instance \cite{ChowTeicher}) provides another justification for considering the average value as a good candidate for a prediction of an observable quantity, since the arithmetic mean is distributed by a Gaussian. Then, arithmetic mean is expected to be observed in the maximum of likely, which is the averaged value.

 When physical observables are elements of a convex space of operators, one can consider arbitrary convex combinations. Then for any two operators $A_1(x,y)$ and $A_2(x,y)$ with the same average value $B(x)$, any combination of the form $B(x,y)=\lambda_1 A_1(x,y) \,+\lambda_2 A_2(x,y)$ with $\lambda_1+\lambda_2=1$ has also as average $B(x)$. This is an invariant property or symmetric property of the average $B(x)$.

 There is plenty of applications of the notion and properties of average. In this note, we discuss two examples of averaging procedures directly related with Finsler geometry: an averaging procedure for geometric operators and some averaging methods for Finsler metrics. We stress the significance of the invariance property of the average discussed above for the case of the averaging in Finsler geometry.

\subsection*{Averaging procedure for a family of linear operators}
Let ${\bf M}$ be an $n$-dimensional manifold and $\pi_{\Sigma}:{\bf \Sigma}\to {\bf M}$ be a sub-bundle of the tangent bundle {\bf TM} such that each of the manifolds  ${\bf \Sigma}_x\hookrightarrow {\bf T}_x{\bf M}$ is compact. We denote by $\pi^*{\bf T}^{(p,q)}{\bf M}$ the tensor bundle over ${\bf M}$ of tensors of type $(p,q)$. The pull-back
bundle $\pi_1 : \pi^*{\bf T}^{(p,q)}{\bf M}\rightarrow{\bf \Sigma}$ is defined by the
commutative diagram
\begin{displaymath}
\xymatrix{\pi^*{\bf T}^{(p,q)}{\bf M} \ar[d]_{\pi_1} \ar[r]^{\pi_2} &
{\bf T}^{(p,q)}{\bf M} \ar[d]^{\pi}\\
{\bf \Sigma} \ar[r]^{\pi_{\Sigma}} & {\bf M}}
\end{displaymath}
Given a family of operators
\begin{displaymath}
A_w:=\{A_w:\pi ^*_w  {\bf T}^{(p,q)}{\bf M}
\longrightarrow \pi^*_w  {\bf 
T}^{(p,q)}{\bf M},\,\,w\in \pi ^{-1}_{{\bf \Sigma}}(x)\,x\in{\bf M}\},
\end{displaymath}
the {\it average value} is defined to be the operator
\begin{displaymath}
\langle A\rangle (x) : {\bf T }^{(p,q)}_x {\bf M} \longrightarrow  {\bf
T}^{(p,q)}_x{\bf M},\quad x\in {\bf M},
\end{displaymath}
determined by the expression
\begin{align}
\langle A_w\,\rangle\,:=\,\frac{1}{vol({\bf \Sigma}_x)}
\Big(\int _{{\bf \Sigma}_x} d\mu (x,y)\pi _2 |_{(x,y)}  
A_{(x,y)}  \pi^* _{(x,y)}\,  \Big)S_x,\quad S_x\in \,{\bf T}^{(p,q)}_x{\bf M}.
\end{align}
Note that this expression implies the existence of a fiber-valued measure on the fiber $\pi^{-1}(x)$.
The volume function is defined as
\begin{align}
vol({\bf \Sigma}_x):=\int_{{\bf \Sigma}_x} 1d\mu(x,y) ,\,\,\,u\in
{\pi^{-1}(x)},\, S_x \subset {\bf T}^{(p,q)}_x {\bf M}.
\end{align}
It requires also the existence of a measure $d\mu(x,y)$ in ${\bf \Sigma}_x$.

In the case of a Finsler structure $({\bf M},F)$, there are several
measures and volume functions on ${\bf T}_x{\bf M}$
that one can define.
We will use the Riemannian measure and volume function induced from
the volume form of the Riemannian metric $({\bf T}_x{\bf M}\setminus \{0\},g_x )$.
This Riemannian metric is defined in local coordinates on ${\bf T}_x{\bf M}\setminus \{0\}$ by the expression
\begin{align*}
g_x:=g_{ij}(x,y)dy^i\otimes dy^j,\quad i,j=1,...,n,
\end{align*}
 with fixed $x\in {\bf M}$
and $y\in {\bf T}_x{\bf M}\setminus \{0\}$.
Its volume form is
\begin{align}
\widetilde{d\mu}(x,y)=\,\sqrt{g(x,y)}\,dy^1\wedge\cdot\cdot\cdot\wedge\,dy^n.
\label{standardvolumeform}
\end{align}
For the sub-manifold $e:{\bf \Sigma}_x\hookrightarrow {\bf T}_x{\bf M}$, the volume form on ${\bf \Sigma}_x$ is the induced  from  $\widetilde{d\mu}$. Also, not that one can still use homogeneous coordinates $(x,y)$, but that they will be necessarily constraint when applied on ${\bf \Sigma}$.

\subsection*{ Averaging of Finsler metrics} In Finsler geometry, one can find examples of averaging procedures. Given a Finsler metric $({\bf M},F)$, the indicatrix on $x$ is the manifold ${\bf I}_x:=\{y\in\,{\bf T}_x{\bf M},\,F(x,y)=1\}$. When integrating the fundamental tensor $g_{ij}:=\frac{1}{2}\,\frac{\partial ^2 F^2(x,y)}{\partial y^i\partial y^j}$ on the indicatrix ${\bf I}_x$, one obtains the components of the {\it averaged metric} are given in local coordinates by the expression
    \begin{align}
    h_{ij}(x):=\langle g_{ij}\rangle(x):=\frac{1}{vol({\bf I}_x)}\Big(\int _{{\bf I}_x} d\mu (x,y)\,g_{ij}(x,y)\Big),\quad i,j=1,...,n.
    \label{averagedmetric}
    \end{align}
The measure $d\mu (x,y)$ and volume form are the standard one induced from the volume form \eqref{standardvolumeform} (see for instance \cite{BCS}). The properties of the averaged metric \eqref{averagedmetric} have been investigated in \cite{Gallego05}. It is particularly useful to investigate generalizations of Riemannian results to the Finsler setting.

Another example of averaging in Finsler metrics, in particular for Berwald metrics, is the averaging introduced by Szab\'{o} used in the  classification of Berwald spaces \cite{Szabo}. In this case, one consider the principal bundle whose standard fiber is the group $G$ of linear transformations leaving invariant the indicatrix ${\bf I}_x$. It is known that $G$ must be a compact Lie group \cite{Wang}. Therefore, there is a left-invariant Haar measure $\mu_G$ on the group $G$ and $vol(G)$ using such measure is finite. Let us fix a point $p\in {\bf  M}$. Given a scalar product $\eta:{\bf T}_p{\bf M}\times {\bf T}_p{\bf M}\to {\bf R}$, one defines a new scalar product
\begin{align}
(X,Y):=\int_G\,\eta(\xi(X),\xi(Y))\mu_G(\xi), X,Y\in {\bf T}_p{\bf M}.
\label{szaboproduct}
\end{align}
This scalar product is invariant under the action of the holonomy group $H_p$ of the Berwald connection. Therefore, for ${\bf M}$ path-connected the scalar product \eqref{szaboproduct} can be extended by parallel transport along curves to the whole manifold, defining a Riemannian metric, and the extension does not depend on the curve connecting $p$ to the generic point $x\in {\bf M}$.

Contrary to the averaged metric $h$, Szab\'o construction cannot be extended to arbitrary Finsler spaces, because in that case the group $G$ is not necessarily compact. This is another reason to consider the construction \eqref{averagedmetric} interesting, since it can be applied to any Finsler structure (indeed, it can also be applied to generalized Finsler metrics, see \cite{MironHrimiucShimadaSabau}).

\section{Convex invariance in Finsler geometry}

Let us consider the $n$-forms\footnote{recall that $d\mu(x,y)$ is an $n-1$-form, the induced volume form on ${\bf I}_x$ of the Riemannian volume form $\widetilde{d\mu}$, eq. \eqref{standardvolumeform}} $\{\Omega^i(x,y)=\,\omega^i(x,y)\wedge d\mu(x,y),\,n=1,...,n\}$ associated to the linear connection $1$-forms $\{\omega^i,\,i=1,...n\}$  of a given linear connection on $\pi^*{\bf TM}$. By fiber integration one obtains the $1$-forms $\{\langle \Omega^i(x,y)\rangle,\,i=1,...,n\}$ on ${\bf M}$ with vector values on ${\bf TM}$,
 \begin{align}
 \langle\Omega^i \rangle(Z)(x):=\,\frac{1}{vol({\bf \Sigma}_x)}
\Big(\int _{{\bf \Sigma}_x} d\mu (x,y)\pi _2 |_{(x,y)}  
\omega^i(x,y)  \pi^* _{(x,y)}\,  \Big)\,\cdot Z,\,i=1,...,n.
\label{averageoperator}
 \end{align}
\begin{proposition}
Let the connection $1$-forms $\omega^i$ on the pull-back bundle $\pi^*{\bf TM}$. Then $\pi^*\langle\Omega^i \rangle$ are the $1$-forms of a linear connection on $\pi^*{\bf TM}$.
\end{proposition}
\begin{proof} Consider
the convex sum of linear forms
\begin{displaymath}
t_1w^i_1( x,y_1)+...t_pw^i_p(x,y_p),\quad t_1+...+t_p=1.
 \end{displaymath}
This convex combination defines a connection on {\bf M}, since convex combinations of connection $1$-forms are connection $1$-forms.
Now let us consider the
manifold ${\bf \Sigma}_x\subset \pi^{-1}(x)\subset {\bf N}$ and a
set of connection $1$-forms $\omega^i(x,y)$ on {\bf M} labeled by points on
${\bf \Sigma}$. We consider integration along the fiber, obtained by limiting process of convex sums attached to finite triangulations \cite{Whitney}. Therefore, the integration operation in \eqref{averageoperator} is defined as the limit of a convex sums. Since all the operations involved in the definition of \eqref{averageoperator} are continuous and the limit exists, the limit
$\{\langle\Omega^i \rangle\}$ defines also the linear connection $1$-forms of a linear connection on {\bf M}. Then we only need to pull-back the limit form, obtaining $\pi^*\langle\Omega^i \rangle$.
\end{proof}
In \cite{Gallego05}, for each linear connection $\nabla$ on $\pi^*{\bf TM}$, it was obtained an averaged connection $\langle\nabla\rangle$. The following results relates the corresponding connection $1$-forms,
\begin{theorem}
The $1$-forms  $\{\langle\Omega^i\rangle (x)\}$ are the connections $1$-forms of the averaged connection $\langle\nabla\rangle$. 
\end{theorem}
\begin{proof}
Let us consider the fiber bundle $\pi:{\bf N}\longrightarrow {\bf
M}$ and the pull-back vector bundle $\pi ^{*} {\bf
TM} \longrightarrow {\bf N} $.
Let $\{e_1,...,e_n\}$ be a local frame basis for the sections of
$\Gamma{\bf TM}$. Then
$\{\pi^*e_1,...,\pi^*e_n\}$ is a basis for the fiber
$\pi^{-1}_u\subset \pi^*{\bf TM},\, u\in {\bf N}$;
$\{h_1,...,h_n\}$ is a basis of the horizontal distribution
$\mathcal{H}_u\,\subset {\bf T}_u{\bf N}$, while $\{v_1,...,v_n\}$
is a basis frame for the vertical distribution $\mathcal{V}_u\,\subset {\bf T}_u{\bf N}$.
Given a non-linear connection $\mathcal{N}$ on
${\bf TN}\longrightarrow {\bf N}$, there are several {\it distinguished} linear connections on the
pull-back bundle $\pi^*{\bf TM}\rightarrow {\bf N}$ constructed using $\mathcal{N}$ (see for instance \cite{Bao, BCS, MironHrimiucShimadaSabau}. Let us consider the connection $1$-forms $\{\omega^i,\,i=1,...,n\}$. Since the connection on $\pi^*{\bf TM}$ is linear, one can define the corresponding Christoffel symbols $\Gamma^i\,_{jk}(x,y)$.
To verify this fact, one needs to check that the connection corresponding to the $1$-forms $\langle\Omega^i \rangle(x)$ and the averaged
connection $\langle \nabla\rangle$ have the same coefficients in a given local frame.
 The result for the connection coefficients of both connections are
\begin{equation}
\langle \Gamma^i\,_{jk}\rangle=\frac{1}{vol({\bf \Sigma}_x)}\,\Big(\int_{{\bf \Sigma}_x} \Gamma^i\,_{jk} (x,y) d\mu(x,y)\Big),
\end{equation}
that coincide with the Christoffel symbols of the connection $\langle \nabla \rangle$ (see reference \cite{Gallego05}).
\end{proof}

The following result shows the existence of an homotopy invariance in Finsler geometry,
\begin{theorem}
There exists a homotopy between each linear connection on the bundle $\pi^*{\bf TM}\longrightarrow {\bf I}$ and the pull-back of its averaged connection.
\end{theorem}
\begin{proof}
Let us consider the connection $1$-forms $w^i$. Then for the average connection the following convex homotopy property holds:
\begin{displaymath}
\langle ((1-t)w^i +t\,\pi^* \langle\Omega^i \rangle )\rangle=\langle\Omega^i \rangle\Rightarrow\,
\pi^*\langle((1-t)w^i +t\,\pi^* \langle\Omega^i \rangle )\rangle =\pi^*\langle\Omega^i \rangle.
\end{displaymath}
 Since the parameter $t\in [0,1]$, we have proved that $\phi$ defines a retraction between $\{(1-t)w^i +t\,\pi^*\omega^i \}$ and $\pi^*\langle\Omega^i \rangle$. Therefore there is a retraction between $w^i$ and $\pi^*\langle\Omega^i \rangle$.
\end{proof}
This homotopy invariance can be interpreted as the invariance of the average value of physical observables discussed in the introduction. As in such case, it shows that the average is not the only information contained in the original Finsler metric or structure. However, it opens the possibility  to explore Finsler spaces from a new perspective. In particular, it makes natural the suggestion that Finsler Geometry is characterized by the following data,
\begin{itemize}
\item {\it A set of associate affine connections}. They are averaged connections, obtained averaging the
Chern, Cartan and Berwald connections. Note that the average of the Chern, Cartan and Berwald connection are in principle different.

\item {\it Riemannian properties}. The averaging of the fundamental tensor of a Finsler structure provides a Riemannian structure. A property is Riemannian if it holds for the full class of Finsler metrics with the same Riemannian metric as average metric.

\item {\it Non-convex invariant properties}. In particular
non-reversibility of the metric function. These properties are lost after averaging.
\end{itemize}

The non-Riemannian and non-affine properties are encoded in two tensors.
The first tensor is the difference between the Finsler metric and the averaged Riemannian structure,
\begin{align}
\delta g:=(g_{ij}-\langle g_{ij}\rangle )\,\pi^* dx^i\otimes \pi^*dx^j.
\label{diferenceofthemetrics}
\end{align}
The second tensor depends on the difference between the initial and averaged connection and measures the degree of non-affine geometry of the Finsler geometry,
\begin{equation}
\delta \Gamma:=\,(\Gamma^i\,_{jk}- \langle \Gamma^i\,_{jk}\rangle)\pi^*\partial_i\otimes \pi^*dx^j \otimes \pi^*dx^k.
\end{equation}
There is a third tensor which measures how much different are the average of the metric and the average of the connection. Let us consider the Levi-Civita connection $^h\nabla$ of the Riemannian metric $h=\langle g_{ij}\rangle\, dx^i\otimes dx^j$ and the corresponding Christoffel symbols $^h\nabla^i_{jk}$. Then one can define the tensor on ${\bf M}$ whose components are given by the torsion
\begin{equation}
D:= \,(^h\Gamma^i\,_{jk}- \langle \Gamma^i\,_{jk}\rangle)\pi^*\partial_i\otimes \pi^*dx^j \otimes \pi^*dx^k.
\end{equation}
This last tensor lives directly on the manifold {\bf M}. Therefore it is a convex invariant quantity.

There are some direct consequences about the above tensors:
\begin{proposition}
The following relations hold,
\begin{align}
\pi^*\langle D\rangle=D,\quad \langle\delta g\rangle=0, \quad \langle\delta \Gamma\rangle=0,
\end{align}
\end{proposition}
\begin{proof}
The first property is proved after a short calculation,
\begin{align*}
\pi^* \langle D\rangle &=\,\pi^*\Big((^h\Gamma^i\,_{jk}- \langle \Gamma^i\,_{jk}\rangle)\partial_i\otimes dx^j \otimes dx^k\Big)\\
& =\,(^h\Gamma^i\,_{jk}- \langle \Gamma^i\,_{jk}\rangle)\pi^*\partial_i\otimes \pi^*dx^j \otimes \pi^*dx^k=\,D.
\end{align*}
The second and third properties are proved analogously.
\end{proof}
A Berwald space is a Finsler space such that the Chern connection defines an affine connection directly on {\bf M}. Then the following result is direct,
\begin{proposition}
Let $({\bf M}, F)$ a Finsler structure.
The following holds,
\begin{enumerate}
\item $\delta \Gamma =0$ iff the structure $({\bf M}, F)$ is Berwald \cite{RF}.

\item $\delta g =0$ iff the structure $({\bf M}, F)$ is Riemannian.
\end{enumerate}
\end{proposition}
 This result suggests that {\it  generalized Finsler spaces} are  characterized by the set of tensors $(\langle g\rangle, \langle \Gamma\rangle,\, \delta g, \delta \Gamma, T)$ such that $\langle g \rangle$ is a Riemannian metric, $\langle \Gamma \rangle$ is an affine connection, $\delta g$ a symmetric tensor with norm  less than $1$ and $T^i_{jk}$ a symmetric tensor in $jk$.

It is conceptually useful to classify all the Finsler metrics with the same average as equivalent. Let $({\bf M},\mathcal{F})$ be the space of Finsler functions on {\bf M},
\begin{definition}
The set of Finsler structures $\{({\bf M},F)\}$ such that they have the same average metric $h$ forms an equivalence class $[F]$ in $({\bf M},\mathcal{F})$ by the relation of having the same average metric.
\end{definition}
That this defines an equivalence relation is direct, since it is transitive, symmetric and reflexive relation. Also, note that for two different Riemannian metrics $[h_1]\cup [h_2]=\,\emptyset$.

It is direct that the class of equivalence of  metric forms a convex set in the space of all the Finsler metrics on ${\bf M}$. Therefore, given a Riemannian metric $h$, which is found to be the average of a fundamental tensor $g$, there is a convex neighborhood $U\subset ({\bf M},\mathcal{F})$ that contains $h$.

\end{document}